\documentclass[conference]{IEEEtran}
\IEEEoverridecommandlockouts
\usepackage{cite}
\usepackage{amsmath,amssymb,amsfonts}
\usepackage{algorithmic}
\usepackage{float}
\usepackage{graphicx}
\usepackage{textcomp}
\usepackage{xcolor}
\newtheorem{theorem}{Theorem}[section]
\newtheorem{corollary}{Corollary}[theorem]
\newtheorem{lemma}[theorem]{Lemma}
\newtheorem{proposition}{Proposition}

\def\BibTeX{{\rm B\kern-.05em{\sci\kern-.025em b}\kern-.08em
    T\kern-.1667em\lower.7ex\hbox{E}\kern-.125emX}}

\begin{document}

\title{Optimal Regulators in Geometric Robotics\\
}

\author{\IEEEauthorblockN{Bousclet Anis}
\IEEEauthorblockA{\textit{Control Engineering Department} \\
\textit{National Polytechnic School}\\
\textit{Student in Master of Engineering} \\
\textit{Process Control Laboratory}\\
Algeria, Algiers \\
anis.bousclet@g.enp.edu.dz }
}
\maketitle

\begin{abstract}
The aim of this paper is to give some existence results of optimal control of robotic systems with a Riemannian geometric view, and derive a formulation of the PMP using the intrinsic geometry of the configuration space.\\ Applying this result to some special cases will give the results of avoidance problems on Riemannian manifolds developed by A. Bloch et al.\\ We derive a formulation of the dynamic programming approach and apply it to the quadratic costs and extend the linear quadratic regulator to robotic systems on Riemannian manifolds and giving an equivalent Riccati equation. We give an optimisation aspect of the Riemannian tracking regulator of F. Bullo and R.M. Murray.\\
Finally, wee apply the theoretical developments to the regulation and tracking of a rigid body attitude. 
\end{abstract}


\section{Introduction}
Robotic manipulators are powerful tools to gain energy, time and money. The problem with the control of this kind of systems is the non-linearity of the dynamics. The classical approach that models the configuration space as an euclidean space and apply the Hamilton's principle to derive the equations and give some regulation results using Lyapunov stability theorems [11] [12] [13] [14] [51]. \newline  \newline On the other hand the work of [8] shows the efficacy of Riemannian geometry in robotics. In [18] [20] [21] [26] [27] the authors propose a Riemannian PD regulator that ensures the regulation of the configuration into a reference position with zero velocity using Lyapunov techniques [3] [10] [12] [13], but the gains can be chosen arbitrarily and are not unique. For linear systems, we solve this problem by fixing a quadratic cost and apply the LQR theory to give the unique optimal control that ensure the regulation. \newline The optimal control is a linear feedback of the state [31] [32] [42], and the proportional coefficient can be computed from an algebraic Riccati equation. \newline Computations in [33] [34] prove that we can recover the Riccati equation by applying the PMP or HJB theory for linear systems with quadratic costs. Work in [43] [44] give a geometric formulation of PMP and HJB theory for state space systems on manifold using the symplectic structure of cotangent bundle of the state space.
\newpage
A. Sacoon et al showed in [39] that the form of the optimal regulator of the rigid body kinematic $R'=R\Omega$ with $R\in SO(3)$ and $\Omega\in \mathfrak{so(3)}=A_{3}(\mathbb{R})$ for the euclidean cost 
$$J(\Omega)=\int_{0}^{\infty}tr(I-R_{d}^{T}.R)+\frac{1}{2}||\Omega||^{2}dt,$$
where $||.||$ is the Frobenius norm,  $$\Omega^{*}(R)=-\frac{R_{d}^{T}R-R^{T}R_{d}}{\sqrt{1+tr(R_{d}^{T}R)}}.$$ We see that this form is not a linear feedback of configuration, this remark prevents to establish an LQR theory for robotic systems using the Euclidean formulation.\\
S. Berkane et al showed in [38] that the optimal control of the Riemannian cost $$J(\Omega)=\int_{0}^{\infty}\frac{1}{2}||\log(R_{d}^{T}R)||^{2}+\frac{1}{2}||\Omega||^{2}dt,$$
is $$\Omega(R)=-k\log(R),$$
for some $k>0$ solution of an algebraic Riccati equation, This result encourages us to use the geometric formulation for optimal control of robotic systems. \\
Works in [40] of C. Liu et al consists to prove that for rigid body dynamics 
$$R'=R\Omega,$$
$$\omega'=I^{-1}(I\omega\times\omega)+\tau.$$
The optimal control for the riemannian cost 
$$J(u)=\int_{0}^{\infty}\frac{1}{2}|\log(R_{d}^{T}R)|_{I}^{2}+\frac{1}{2}|\Omega|_{I}^{2}+\frac{1}{2}|[\tau]_{\times}|_{I}^{2}dt$$
is $\tau=-k_{P}\log(R_{d}^{T}R)-k_{D}\omega$ with $k_{P},k_{D}$ are solution of algebraic Riccati equation. \\ 
The goal here is to prove that we can establish a geometric LQR theory for robotic systems, by proving that for riemannian cost, the optimal control is a riemannian PD regulator. 
\section{Mathematical Preliminaries}
Here we recall some preliminaries about geometric robotic and some functional analysis.
\begin{figure}[H]
  \centering 
      \includegraphics[scale=0.4]{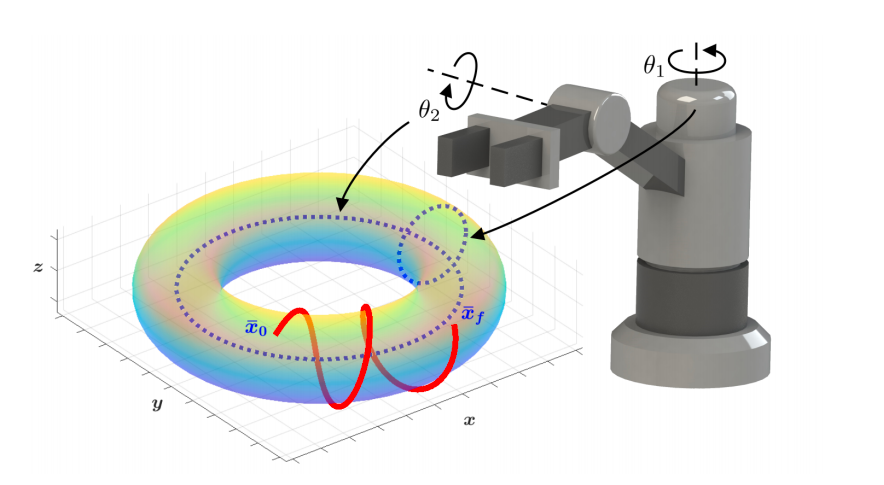} 
      \caption{Configuration space of a robot [52]}

           \label{fig:my_label}
            \end{figure}
The configuration space of a robot is a $n$ dimensional manifold $M\subset (SO_{3}(\mathbb{R})\times\mathbb{R}^{3})^{s}$ where $s$ is the number of rigid bodies that constitute the robot, and $n=6s-r$ is the degree of freedom and $r$ is the number of holonomic constraints.\\
The Riemannian manifold $(M,g)$ with $g=<,>$ is the kinetic energy of the robot allow us to express trajectories of the robot subject to the potential $W\in C^{\infty}(M)$ and the control force $u\in L^{1}_{loc}(I,T^{*}M)$ such that $u(t)\in T_{\gamma(t)}M^{*}$ by $$\frac{D\dot{\gamma}}{Dt}=-grad_{g}(W)(\gamma)+g_{\gamma}^{\#}(u)$$ where $\frac{D}{Dt}:\Gamma(\gamma)\rightarrow{\Gamma(\gamma)}$ is the covariant derivative.\\ $grad_{g}(W)$ is the covariant gradient i.e the vector field that ensure $dW_{q}(v)=<grad_{g}(W)(q),v>$ and $g_{q}^{\#}$ is the tangent cotangent isomorphism.\\
We call a state the element $(q,v)\in TM$, that is a state is the combination of a configuration and a velocity. \\
The curvature tensor is $R:\Gamma(TM)^{3}\rightarrow{\Gamma(TM)}$ defined by $$R(X,Y)Z=D_{X}D_{Y}Z-D_{Y}D_{X}Z-D_{[X,Y]}Z$$ with $D:\Gamma(TM)\rightarrow{\Gamma(TM)}$ is the Levi-Civita connexion.\\
A locally absolutely continuous function $f\in LAC(I,\mathbb{R})$ is a function such that $f(t)=f(t_{0})+\int_{t_{0}}^{t}g(s)ds$ for $g\in L^{1}_{loc}(I)$, from the Lebesgue theorem [25], an absolutely continuous function is a.e differentiable and $f'=g$ a.e, we have this result [42] [33]
\begin{theorem}{Cauchy-Lipschitz in Optimal Control}\\
for $u\in \Gamma^{1}_{loc}(\gamma)$, the dynamic $\frac{D\dot{\gamma}}{Dt}=-grad(U)(\gamma)+u$ have a unique LAD local solution for each initial conditions on $TM$ defined in an open interval.
\end{theorem} 
From Mazur's and the converse of Lebesgue dominated convergence theorems, we have the following facts, if $\Omega\subset\mathbb{R}^{m}$ convex and closed, $L^{2}([0,T],\Omega)$ is weakly closed, and if $\Omega$ is convex and compact, $L^{2}([0,T],\Omega)$ is weakly compact [29].\\\\
We denote by $U\in C(M\times M)$ by the minimal kinetic energy that takes the robot from a configuration to another, namely $U(p,q)=\frac{1}{2}d_{g}^{2}(p,q)$.\\
For $q^{*}\in M$ the function $V(.)=U(.,q^{*})$ is smooth on a normal neighborhood of $q^{*}$, and its gradient is $\nabla V(q)=-\exp_{q}^{-1}(q^{*})$, in fact we have $U(q^{*},q)=\frac{1}{2}|\exp_{q}^{-1}(q^{*})|_{g}^{2}$.\\
We denote $inj_{q}M$ by the injectivity radius of $(M,g)$ at $q$ by the biggest number $r>0$ such that $\exp_{q}:B(0,r)\rightarrow{M}$ is a diffeomorphism.\\
If $d_{g}(p,q)<inj_{q}M$ there is a unique minimizing geodesic that joints $p$ to $q$.
\begin{figure}[H]
  \centering 
      \includegraphics[scale=0.15]{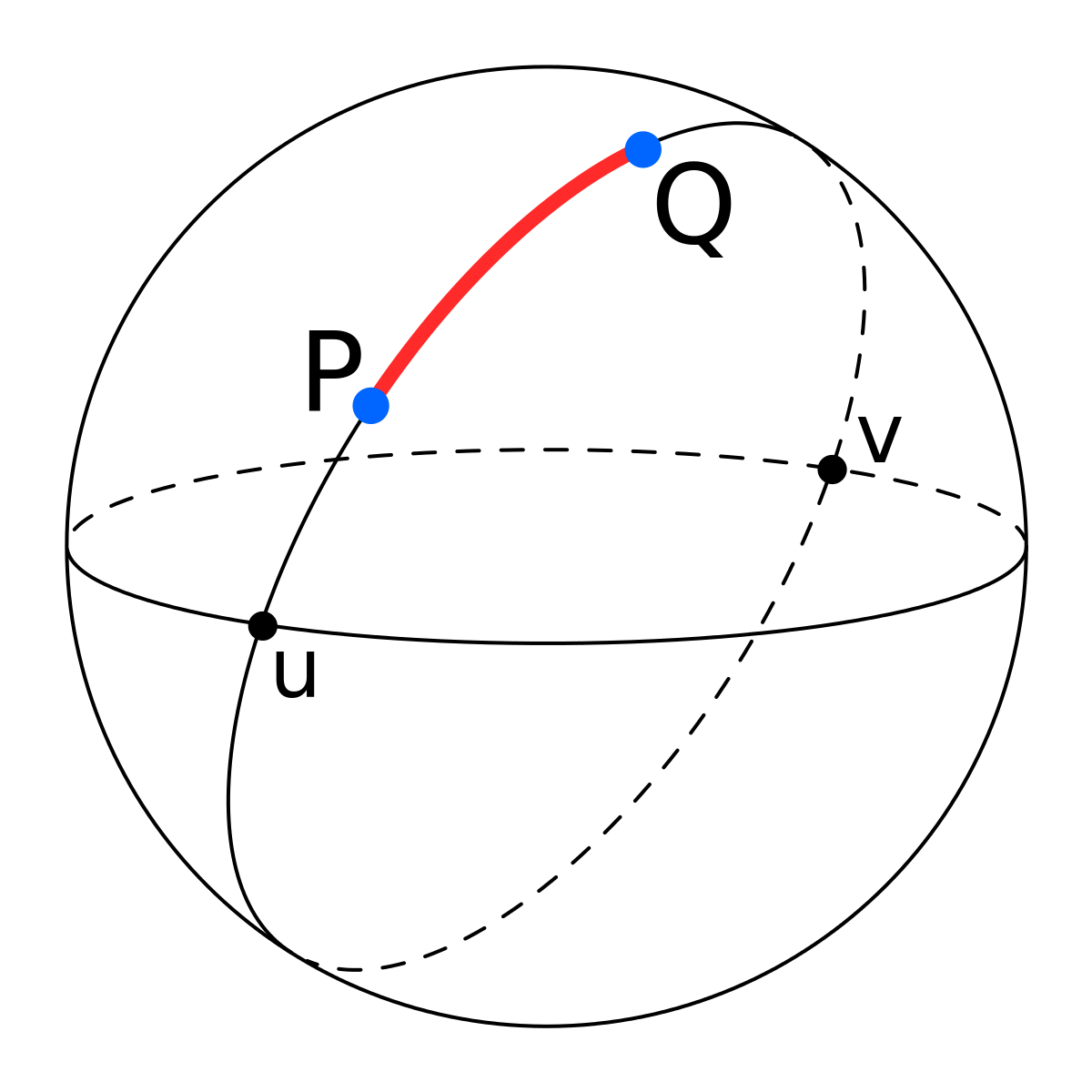} 
      \caption{Geodesic distance}

           \label{fig:my_label}
            \end{figure}
We will show that for natural cost, the optimal control is a Riemannian PD feedback, where the proportional action is exactly $\exp^{-1}_{q}(q^{*})$, the major problem of the geometric approach is fastidious computations, in practice it is very difficult to compute the exponential map neither its inverse, so we do two approximations \\\\
1)- the case where $M$ is an arbitrary manifold $$\exp^{-1}_{q}(q^{*})=q^{*}-q$$. \\
2) the case where $M=SO(3)$, we take $$\exp^{-1}_{R}(R_{d})=R.\log(R^{T}.R_{d})$$ where $\log :\left\{R\in SO(3), tr(R)\neq -1\right\}\rightarrow{A_{3}(\mathbb{R})}$ is the logarithm that can be computed easily from Rodriguez formula [49] [21].\\
$$\log(R)=\frac{\phi(R)}{\sin(\phi(R))}Skew(R)$$ with $Skew(R)=\frac{R-R^{T}}{2}$ is the projection of $R$ on the space of anti-symmetric matrices, and $\phi(R)=\arccos(\frac{tr(R)-1}{2})$.\\\\
The motivation of the approximation on $SO(3)$ is that for left invariant metric, knowledge of the exponential map at $I_{3}$ gives us its values in all $TSO(3)$. On the other hand, when the metric is bi invariant, the exponential map of the metric and the Lie exponential (matrix exponential) are the same, this allows us to approximate the two logarithms. \\
We see that this allow us to recover the classical LQR theory, and classical results of Linear regulation, or results about regulation of rigid bodies. \\
It's clear that $L(q,v)=\frac{1}{2}d_{g}^{2}(q^{*},q)+\frac{1}{2}|v|^{2}$ is a Lyapunov function for the dynamic $$\frac{Dq'}{Dt}=\exp_{q}^{-1}(q^{*})-k.v$$ where $k>0$ and $(q_{0},v_{0})$ satisfies $L(q_{0},v_{0})\leq \frac{1}{2}inj_{q^{*}}(M)^{2}$.\\ \\ 
We will prove that it's an optimal feedback regulator for a natural cost.
\begin{figure}[H]
  \centering 
      \includegraphics[scale=0.5]{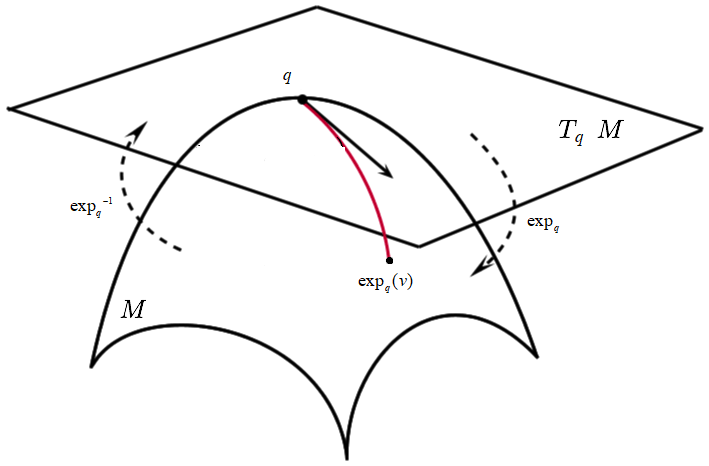} 
      \caption{Riemannian logarithm}

           \label{fig:my_label}
            \end{figure}
\section{Existence Results in Optimal Control of Robotic Systems}
The aim of this part is to give some existence results of optimal control laws by functional analysis tools [29], we adapt the results in [31] [33] [41] [42] to the formulation in geometric robotics.\\
we need the following result
\begin{lemma}{Inverse of the covariante derivative}
Let $\gamma:[0,T]\rightarrow{M}$ be LAD, so there exists a unique linear operator $I_{\gamma}:L^{1}(\gamma)\rightarrow{AC(\gamma)}$ that satisfies $I(\frac{DV}{Dt})(t)=V(t)-P_{0}^{t}(V(0))$ for all $t\in I$ and all $V\in AC(\gamma)$, and $\frac{DI(w)}{Dt}(t)=W(t)$ a.e in $I$ for all $W\in L^{1}(\gamma)$.
\end{lemma}
\textbf{Proof :}\\
We construct this operator locally, and we extend it by the uniqueness, we solve the equation
$$\dot{W}^{k}+\sum_{i,j=1}^{n}\Gamma_{i,j}^{k}(\gamma)\dot{\gamma}^{i}W^{j}=V^{k}$$ with initial condition $W(0)=0$, the linear Cauchy-Lipschitz theorem gives us this solution that satisfies $\frac{DI(V)}{Dt}=V$, we see also that $W=V-P_{0}^{t}(V(0))$ is a solution of this problem, so the conclusion follows.
\subsection{Time-optimal control result}
We denote $U_{T}\subset L^{2}([0,T],\Omega)$ by the set of controls such that the robot initialized in $(q_{0},v_{0})\in TM$ the trajectory is defined in $[0,T]$, and the trajectories lies in a fixed compact $K\subset TM$. \\
let $E:U_{T}\rightarrow{TM}$ defined by $E(u)=(q_{u}(T),q_{u}'(T))$, the accessible set is $E(U_{T})$.\\
The fact that trajectories are lies in a fixed compact $K$ and the control lies in a convex compact implies that the accessible set is compact and depend continuously for the Hausdorff distance.
\fcolorbox{black}{lightgray}{
\begin{minipage}{0.48\textwidth}

\begin{theorem}
Suppose that $\Omega$ is convex and compact, and consider the dynamic $$\frac{Dq'}{Dt}=-grad(W)(q)+\sum_{j=1}^{m}u_{j}f_{j}(q)$$ so the accessible set $Acc_{t}(q_{0},v_{0})$ is compact and depend continuously for the Hausdorff distance.
\end{theorem}
\end{minipage}}\\ \\
\textbf{Proof :}\\
We are interested by the dynamic 
 $$q'=v$$ $$\frac{Dv}{Dt}=-\nabla W(q)+\sum_{j=1}^{m}f_{j}(q)u_{j}.$$
 We prove the compactness of $Acc_{t}$, let $(q_{n}(t),v_{n}(t))$ be a sequence of $Acc_{t}$, because $u_{n}\in L^{2}([0,t],\Omega)$ and that $\Omega$ is bounded. $||u_{n}||_{L^{2}}$ is bounded, then, there exist $u\in L^{2}([0,t],\Omega)$ such that $u_{\phi(n)}$ converge weakly to $u$.\\
 On the other hand $(q_{\phi(n)},v_{\phi(n)})$ is bounded in $H^{1}$, so it converges weakly for some $(q,v):[0,t]\rightarrow{TM}$, we know that $H^{1}([0,t])$ can be compactly injected in $C^{0}([0,t])$, this proves the uniform convergence. \\
 So we have 
 $$ 
     q(t)=q_{0}+\int_{0}^{t}v(s)ds, 
$$
$$
     v(t)=P_{0}^{t}(v_{0})+I_{q}(-\nabla W(q)+\sum_{j=1}^{m}f_{j}(q)u_{j})(t).
$$
 
 This conclude that $u$ is the control that gives $(q,v)$.\\
 Now we prove the continuous dependence\\
 Let $t,s\in [0,T]$, $\epsilon >0$. it suffices to prove that there exists $\delta>0$ such that if $|t-s|\leq\delta$, so for $y\in Acc_{s}$ there exists $x\in Acc_{t}$ such that $||x-y||\leq \epsilon$.
  We fix $t<s$, for the $y\in Acc_{s}$, we choose $x=(q_{u}(t),v_{u}(t))$, we have
  $$y-x=(\int_{t}^{s}v(\tau)d\tau,[P_{0}(v_{0})+I_{q}(-\nabla W+f_{j}u_{j})]|_{t}^{s}).$$
 The fact that quantities that are in the integral are locally integrable concludes. \\ 
 \textbf{Remarque :} The hypothesis about the trajectories that must be in a fixed compact is natural and essential in the proof.\\
We have this corollary \\
\fcolorbox{black}{lightgray}{
\begin{minipage}{0.48\textwidth}

\begin{corollary}
If $(q_{1},v_{1})\in Acc_{T}$ for some $T>0$, so there exist a time-optimal control that takes $(q_{0},v_{0})$ to $(q_{1},v_{1})$ in time $t^{*}$, in addition $(q_{1},v_{1})\in \partial Acc_{t^{*}}(q_{0},v_{0})$. 
\end{corollary}
\end{minipage}}\\ \\
\textbf{Proof :}\\
 Let $t^{*}=\inf\left\{t>0, (q_{1},v_{1})\in Acc_{t}(q_{0},v_{0})\right\}$, suppose that $(q_{1},v_{1})\notin Acc_{t^{*}}(q_{0},v_{0})$, let $\epsilon=\frac{1}{2}d(x_{1},x_{1}')$ (where $x_{1}'$ satisfies the minimal distance, this point exists because the accessible sets are compacts). so there exists $\delta>0$ such that if $t\in[t^{*},t^{*}+\delta]$ for every $x\in Acc_{t}$, there exists $y\in Acc_{t^{*}}$ such that $||x-y||\leq\epsilon$. From the definition of infimum , there exists $t^{*}<t<t^{*}+\delta$ that satisfies $x_{1}\in Acc_{t}$, we take $x=x_{1}$ the corresponding $y$ is in $Acc_{t}$ closer from $x_{1}$ than $x_{1}'$, this is a contradiction.\\\\
 If $x_{1}\notin \partial Acc_{t^{*}}$, so $x_{1}\in \mathring{Acc_{t^{*}}}$ and by continuity there exists $t<t^{*}$ such that $x_{1}\in Acc_{t}$, this contradicts the time optimality. 
 \begin{figure}[H]
  \centering 
      \includegraphics[scale=0.7]{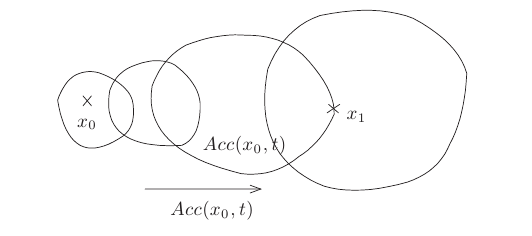} 
      \caption{Temps optimalité}

           \label{fig:my_label}
            \end{figure}
\subsection{Cost optimal control results}
Now we give existence results for optimisation of some cost, we note $U=\frac{1}{2}d_{g}^{2}$ and $P$ is the parallel transportation along the unique minimizing geodesic (this is well defined for sufficiently close configurations), let $R\in S_{n}^{++}(\mathbb{R})$. 
\begin{theorem}
Consider the dynamic $\frac{Dq'}{Dt}=-grad(W)(q)+\sum_{j=1}^{m}u_{j}f_{j}(q)$ and $\Omega$ is closed and convex, so there exist an optimal control for the cost $$J(u)=\int_{0}^{T}||u||_{R}^{2}+\frac{1}{2}|v_{u}|^{2}+U(q_{1},q_{u})dt$$ 
\end{theorem}

\textbf{Proof :}\\
Let $\beta>0$ such that $\beta\leq R$, so $\beta.||u||^{2}\leq J(u)$, let $u_{n}\in U_{T}$ be a minimizing sequence i.e $J(u_{n})\rightarrow{\inf_{u\in U_{T}} J(u)}$, this sequence is clearly bounded, so there exists a subsequence that converges to $u\in L^{2}([0,T],\Omega)$. We prove now that $u\in U_{T}$, we recall that $L^{2}([0,T],\Omega)$ is weakly closed.\\
Let $(q_{n},v_{n})$ be the trajectories with control $u_{n}$ initialized at $(q_{0},v_{0})$, $$q_{n}(t)=q_{0}+\int_{0}^{t}v_{n}(s)ds,$$ $$v_{n}(t)=P_{0}^{t}(v_{0})+I_{q_{n}}(-\nabla W(q_{n})+f_{j}(q_{n})u_{j,n})(t).$$ 
 The sequence $(q_{n},v_{n})\in H^{1}([0,T])$ is bounded, we can extract a sequence that converges weakly to $(q,v)\in H^{1}([0,T])$, by the Sobolev compactness theorem, this sequence converges uniformly to $(q,v)$ in $C^{0}$, so $$q(t)=q_{0}+\int_{0}^{t}v(s)ds$$ $$v(t)=P_{0}^{t}(v_{0})+I_{q}(-\nabla W(q)+f_{j}(q)u_{j})(t)$$
This confirms that $u\in U_{T}$ because $q,v$ are defined in $[0,T]$, on the other hand, we have $||u||_{R,L^{2}}\leq \liminf ||u_{n}||_{R,L^{2}}$ (the norm $||.||_{R,L^{2}}$ is convex $L^{2}$ and continuous, so it is lower semi-continuous for the weak topology) and thus by uniform convergence of  $(q_{n},v_{n})$ we have $$\inf_{u\in U}J(u)\geq J(u).$$
As for the regulation problem, we can consider a tracking regulator as in this proposition

\begin{proposition}
Consider the dynamic $\frac{Dq'}{Dt}=-grad(W)(q)+\sum_{j=1}^{m}u_{j}f_{j}(q)$ and $\Omega$ is closed and convex, so there exist an optimal control for the cost $$J(u)=\int_{0}^{T}||u||^{2}+\frac{1}{2}|v_{u}-P(q_{u},q_{ref})v_{ref}|^{2}+U(q_{ref},q_{u})dt$$ 
\end{proposition}

\section{Riemannian Formulation of PMP}
Works in [43] [44] give a symplectic formulation of PMP for state space systems on manifolds, here we give a Riemannian formulation of PMP for robotic systems that explicitly implies the geometry of configuration space. \\
we start by giving an estimation of how perturbing a control changes the configuration, and it's here where the Riemannian tensor will appears ! after that we will apply this results to the Riemannian double integrator, and recover results in [45] [46] [47] of Bloch, Silva and Colombo.\\
\subsection{Statement of the problem}
We are interesting by the following optimisation problem\\
$$\frac{Dq'}{Dt}=-grad(W)(q)+\sum_{j=1}^{m}f_{j}(q)u_{j}$$ with $(q(0),q'(0))=(q_{0},v_{0})$ and $u\in U_{T}$, we take $L:TM\times\mathbb{R}^{m}\rightarrow{\mathbb{R}}$ and $g:TM\rightarrow{\mathbb{R}}$ are smooth, and we want to minimize $J:U_{T}\rightarrow{\mathbb{R}}$ $$J(u)=\int_{0}^{T}L(q_{u}(t),v_{u}(t),u(t))dt+g(q_{u}(T),v_{u}(T)).$$
Our problem is to give a necessary condition for $u^{*}\in U_{T}$ that verifies $J(u^{*})=\inf_{u\in U_{T}}J(u)$.
\subsection{Simple variations of control}
The tangent bundle $TM$ of the configuration space (state space) is furnished with the Sasaki metric $G$ [24] [25], for $(v_{0}',v_{1}')\in T_{(q_{0},v_{0})}TM$ $\exp_{(q_{0},v_{0})}(v_{0}',v_{1}')$ is the exponential map.\\
let $s>0$ and $a\in\Omega$, we denote $(q^{*},u^{*})$ by the optimal trajectory, and we define the following perturbation of the optimal control $u_{\epsilon}(t)=a$ if $t\in ]s-\epsilon,s[$ and $u_{\epsilon}(t)=u^{*}(t)$ otherwise, it is clear that $u_{\epsilon}\in L^{2}([0,T],\Omega)$, and we will show that in fact, $u_{\epsilon}\in U_{T}$ for $\epsilon>0$ sufficiently small.\\
the first step is to give the configuration by mean of the optimal configuration when the initial conditions are close.\\
\begin{lemma}{Perturbation analysis}\\
let the dynamic $$\frac{Dq_{\epsilon}'}{Dt}=-grad(W)(q_{\epsilon})+\sum_{j=1}^{m}f_{j}(q_{\epsilon})u_{j}$$ and $(q_{\epsilon},v_{\epsilon})(0)=\exp_{(q_{0},v_{0})}(\epsilon.(v_{0}',v_{1}')+o(\epsilon))$, so $q_{\epsilon}(t)=\exp_{q(t)}(\epsilon Y(t)+o(\epsilon))$ with $Y\in AD(q)$ is the solution of the linear equation $$\frac{D^{2}Y}{Dt^{2}}=-D_{Y}grad(W)(q)+R(q',Y)q'+\sum_{j=1}^{m}D_{Y}f_{j}(q)u_{j}$$
and $Y(0)=v_{0}'$ $\frac{DY}{Dt}(0)=v_{1}'$. 
\end{lemma}
\textbf{Proof :}\\
By continuous dependence on initial conditions, we can write $q_{\epsilon}(t)=\exp_{q(t)}(\epsilon Y(t)+o(\epsilon))$ for some $Y\in AD(q)$. \\
Applying $\frac{D}{D\epsilon}$ on the family of equations that satisfies $q_{\epsilon}$ and using the fact that if we invert two covariant derivative the curvature will appears, we conclude by putting $\epsilon=0$. \\  \\
Now we can estimate the configuration of the robot when there is small variations of the control with the following result :
\begin{lemma}
for the dynamic $$\frac{Dq_{\epsilon}'}{Dt}=-grad(W)(q_{\epsilon})+\sum_{j=1}^{m}f_{j}(q_{\epsilon})u_{\epsilon,j}$$ with $(q_{0},v_{0})$ as initial conditions.\\
so $q_{\epsilon}(t)=\exp_{q(t)}(\epsilon.Y(t)+o(\epsilon))$ with $Y(t)=0$ in $[0,s]$ and $Y$ is the solution of the linear equation $$\frac{D^{2}Y}{Dt^{2}}=-D_{Y}grad(W)(q)+R(q',Y)q'+\sum_{j=1}^{m}D_{Y}f_{j}(q)u_{j}$$ and $Y(s)=0$, $\frac{DY}{Dt}(s)=\sum_{j=1}^{m}f_{j}(q(s)).(a_{j}-u_{j}(s))$
\end{lemma}
\textbf{Proof :}\\
It's clear that $q_{\epsilon}=q$ when $t\in[0,s-\epsilon]$, we easily establish that $(q_{\epsilon},v_{\epsilon})(s)=\exp_{(q_{0},v_{0})}(\epsilon(0,\sum_{j=1}^{m}f_{j}(q(s)).(a_{j}-u_{j}(s)))+o(\epsilon))$. \\
The previous lemma concludes. 
\subsection{Riemannian geometric PMP}
Now we can give an intrinsic formulation of PMP for robotic systems.\\
\fcolorbox{black}{lightgray}{
\begin{minipage}{0.48\textwidth}

\begin{theorem}{Intrinsic formulation of PMP}
consider the robot dynamic $\frac{Dq'}{Dt}=-grad(W)(q)+\sum_{j=1}^{m}f_{j}(q)u_{j}$ with $u\in U_{T}$, so if $u^{*}\in U_{T}$ is optimal for the cost as defined in the statement problem, there exists $p_{1},p_{2}\in AC(q)$ that satisfies Hamilton's ODE $$\frac{Dp_{1}}{Dt}=-A^{*}(p_{2})-grad_{q}L$$ $$\frac{Dp_{2}}{Dt}=-p_{1}-grad_{v}L$$ and $p(T)=(p_{1},p_{2})(T)=\nabla g(q(T),v(T))$, and satisfies the minimization principle a.e $$H(q(t),v(t),p(t),u(t))=\min_{a\in\Omega}H(q(t),v(t),p(t),a)$$ where $H=<p_{1},v>+<p_{2},-grad(W)+f_{j}u_{j}>+L$, and $H$ is constant along these trajectories. 
\end{theorem}
\end{minipage}}\\ \\
\textbf{Proof :}\\
We suppose in first time that $L=0$, so it's a terminal value optimisation problem, $J(u)=g(q_{u}(T),v_{u}(T))$, the following step is to estimate $\frac{d}{d\epsilon}J(u_{\epsilon})|_{\epsilon=0}$\\
using previous lemma, we have $\frac{d}{d\epsilon}J(u_{\epsilon})=dg_{(q,v)(T)}(Y(T),\frac{DY}{Dt}(T))$, define $(p_{1},p_{2})\in AC(q)$ solution of $$
    \frac{Dp_{1}}{Dt}=-A^{*}(p_{2}),
$$ $$
    \frac{Dp_{2}}{Dt}=-p_{1},
$$ with $(p_{1},p_{2})(T)=\nabla g(q(T),v(T))$. so the expression $<p_{1},Y>+<p_{2},\frac{DY}{Dt}>$ is constant. It's derivative is  $$<\frac{Dp_{1}}{Dt},Y>+<p_{1},\frac{DY}{Dt}>+<\frac{Dp_{2}}{Dt},\frac{DY}{Dt}>$$ $$+<p_{2},\frac{D^{2}Y}{Dt^{2}}>=-<A^{*}(p_{2}),Y>+<p_{2},\frac{D^{2}Y}{Dt^{2}}>.$$ Where $A$ is the operator such that
$$\frac{D^{2}Y}{Dt^{2}}=A(Y),$$
we conclude by taking $A^{*}$ the adjoint operator.\\
On a $\frac{d}{d\epsilon}J(u_{\epsilon})|_{\epsilon=0}=\sum_{j=1}^{m}(a_{j}-u_{j}^{*})<p_{2}(s),f_{j}(q(s))>$, the fact that $u^{*}$ is optimal, and $s, a$ are arbitrary, we conclude that $$\sum_{j=1}^{m}u_{j}^{*}<p_{2}(s),f_{j}(q(s))>=\min_{a\in\Omega}\sum_{j=1}^{m}a_{j}<p_{2}(s),f_{j}(q(s))>.$$
\subsection{Avoidance problem in robotics}
Here we apply the developed result to the particular case $$\frac{Dq'}{Dt}=u$$ we want to minimize the cost $$J(u)=\int_{0}^{T}[U(q^{*},q_{u}(t))+\frac{1}{2}|v_{u}(t)|^{2}+\frac{\alpha}{2}|u|^{2}+V(q_{u}(t))]dt$$ with $V:M\rightarrow{\mathbb{R}}$ is the avoidance function, it looks like $V(q)=\sum_{i=1}^{o}\frac{1}{O_{i}(q)}$ with $o$ is the number of obstacles, and the obstacles are $O_{i}=\left\{q\in M, O_{i}(q)\leq 0\right\}$, we have the following result as in [45] [46] [47].
\fcolorbox{black}{white}{
\begin{minipage}{0.48\textwidth}
\begin{theorem}{Avoidance problem}\\
the optimal control of the Riemannian double integrator for the cost $J$ is the solution of the equation $$\frac{D^{2}u}{Dt^{2}}=R(v,u)v+\frac{1}{\alpha}grad(U+V)(q)+\frac{1}{\alpha}u$$ with $u(T)=0$ and $\frac{Du}{Dt}(T)=0$.
\end{theorem}
\end{minipage}}\\ \\
\textbf{Proof :}\\
The Hamiltonian is a strict convex function on $u$, so minimum is the unique critical point of
$$H=<p_{1},v>+<p_{2},u>+U+\frac{1}{2}|v|^{2}+V+\frac{\alpha}{2}|u|^{2}.$$
This gives $u^{*}=-\frac{1}{\alpha}p_{2}$.\\
On the other hand, computations gives $$\nabla L(q,v)=(\nabla (U+V)(q),v),$$
by the fact $<R(X,Y)Z,W>=<R(Z,W)X,Y>$, we have $A^{*}(p_{2})=R(v,p_{2})v$. \\
This concludes. 
\begin{figure}[H]
  \centering 
      \includegraphics[scale=1]{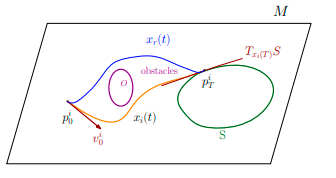} 
      \caption{Avoidance problem [47]}

           \label{fig:my_label}
            \end{figure}
\subsection{Optimal regulation}
we can give a result about regulation in finite time.
\fcolorbox{black}{white}{
\begin{minipage}{0.48\textwidth}

\begin{proposition}{Curvature and Regulation}\\
let the dynamic $\frac{Dq'}{Dt}=u$ and the cost $J(u)=\int_{0}^{T}\frac{\alpha}{2}|u|^{2}+U(q^{*},q_{u}(T))+\frac{1}{2}|v_{u}(T)|^{2}$
then the optimal control is solution of the Jacobi equation $$\frac{D^{2}u}{Dt^{2}}=R(v,u)v$$
with $u(T)=-\frac{1}{\alpha}v(T)$ and $\frac{Du}{Dt}(T)=-\frac{1}{\alpha}grad(U)(q^{*},q(T))$
\end{proposition}
\end{minipage}}\\ \\
\textbf{Proof :}\\
Same computations as for previous theorem. 
\section{HJB theory and applications}
Here we give the formulation of HJB theory in the case of robotic systems, formulation in [43] [44] gives a symplectic one for dynamical systems on manifolds, here we give a Riemannian formulation for robotic systems.
\subsection{Riemannian formulation of HJB}
Let $L:[0,T]\times TM\times\mathbb{R}^{m}\rightarrow{\mathbb{R}}$ and $g:TM\rightarrow{\mathbb{R}}$ be smooth, we are interesting by solving the optimal control problem as exposed in the statement of the problem.\\
as in [9] [38] [37], we can prove that if the value function is smooth $V(t,q,v)=\inf_{u\in U_{T}}\int_{t}^{T}L(s,q_{u},v_{u},u)dt+g(q_{u},v_{u})(T)$ with $q_{u}(t)=q$ and $v_{u}(t)=v$ and $\frac{Dv_{u}}{Dt}=-\nabla W(q_{u})+f_{j}(q_{u})u_{j}$. \\
We introduce the Hamiltonian 
$$H=\inf_{u\in\Omega}\left\{<p_{1},v>+<p_{2},-\nabla W(q)+f_{j}(q)u_{j}>+L\right\}.$$
We have the following theorem.
\fcolorbox{black}{white}{
\begin{minipage}{0.48\textwidth}

\begin{theorem}
We suppose that $V\in C^{1}([0,T]\times TM)$, and the Lagrangian is coercive in $u$ 
$$L\geq \beta |u|^{2}.$$
So it satisfies the Hamilton Jacobi Bellman PDE 
 $$\frac{\partial V}{\partial t}+H(q,v,\nabla V)=0,$$ with terminal value $$V(T,q,v)=g(q,v),$$ with $\nabla V$ is the gradient associated to the Sasaki metric on the tangent bundle. 
\end{theorem}
\end{minipage}}\\ \\
\textbf{Proof :}\\
Let $(q,v)\in TM$ and $t\in[0,T]$, we denote $V:[0,T]\times TM\rightarrow{\mathbb{R}}$ the value function by 
$$
    V=\inf_{u\in U_{\Omega}}[ \int_{t}^{T}L(s,q_{u}(s),v_{u}(s),u(s))ds+g(q_{u}(T),v_{u}(T))], 
$$ where $(q_{u},v_{u})$ are solutions of the robot dynamic, and $(q_{u}(t),v_{u}(t))=(q,v)$.\\
for $t\leq s\leq t+h$ we take $u(s)=a\in\Omega$, and for $s>t+h$ we choose an optimal control, the cost of this control is $$\int_{t}^{t+h}L(s,q(s),v(s),u(s))ds+V(t+h,q(t+h),v(t+h))$$ $$\geq  V(t,q(t),v(t),$$ so we have $$\frac{V(t+h,q(t+h),v(t+h)-V(t,q(t),v(t))}{h}$$ $$+\frac{1}{h}\int_{t}^{t+h}L(s,q(s),v(s),a)ds\geq 0.$$ let $h\rightarrow{0}$
$$\frac{\partial V}{\partial t}+<\nabla_{q} V,q'>+<\nabla_{v}V,\frac{Dv}{Dt}>$$ $$+L\geq 0,$$ and this gives $$\frac{\partial V}{\partial t}+\inf_{a\in\Omega}[<\nabla_{q}V,v>+<\nabla_{v}V,-\nabla W(q)+f_{j}(q)a_{j}>$$ $$+L(t,q,v,a)]\geq 0.$$ By choosing an optimal control we obtain
$$\frac{\partial V}{\partial t}+\min_{a\in\Omega}[<\nabla_{q}V,v>+<\nabla_{v}V,-\nabla W(q)+f_{j}(q)a_{j}>$$ $$+L(t,q,v,a)]=0.$$ 
and this is the Riemannian formulation of HJB equation. \\\\
One can also prove that under some suppositions, $V$ is the unique viscosity solution
[33] [30] of the HJB equation. 
\subsection{Strategy of optimal control by HJB}
The strategy of optimal control by HJB is the following, we compute the minimum of the Hamiltonian $$H=<p_{1},v>+<p_{2},-\nabla W(q)+f_{j}(q)u_{j}>+L(t,q,v,u)$$ and the argument is $u^{*}$. We solve HJB equation and express the optimal control as a feedback.\\\\
We verify now that this construction gives us an optimal control\\ $$\int_{0}^{T}L(s,q(s),v(s),u(s))+g(q(T),v(T))=-\int_{0}^{T}[\frac{\partial V}{\partial t}$$ $$+<\nabla_{q}V,v> +<\nabla_{v}V,-\nabla W(q)+f_{j}(q)u_{j}>]ds$$ $$+g(q(T),v(T)) 
=-\int_{0}^{T}\frac{d}{ds}V(s,q(s),v(s))ds+g(q(T),v(T)$$ $$=V(0,q,v).
$$ This confirms optimality.\\ \\
The optimal control verify then $$\frac{d}{d}V(t,q(t),v(t))=-L(t,q(t),v(t),u(t)),$$
In particular when $L\geq 0$, the value function is non increasing along optimal trajectories. 

\subsection{Optimal Regulation (Riemannian-LQR)}
Here we look for a regulator that takes the configuration to a reference one $q^{*}$ with zero velocity, and minimize the cost $J(u)=\int_{0}^{\infty}[U(q^{*},q_{u}(t))+\frac{1}{2}|v_{u}(t)|^{2}+\frac{\alpha}{2}|u|^{2}]dt$ for the dynamic $\frac{Dv_{u}}{Dt}=u$.\\
we will prove that the optimal control is $u=-k_{P}.\nabla U-k_{D}v$ for some $k_{P},k_{D}>0$ that solves an algebraic Riccati equation.
\fcolorbox{black}{lightgray}{
\begin{minipage}{0.48\textwidth}

\begin{theorem}{LQR-Theory}\\
the optimal control of the Riemannian double integrator for the natural cost is a PD regulator $u(q,v)=-\frac{k_{3}}{\alpha}\nabla U(q)-\frac{k_{2}}{\alpha}v$ with $(\frac{k_{3}}{\alpha},\frac{k_{2}}{\alpha})=R^{-1}.B^{T}.K$ and $K$ is the unique positive definite matrix that satisfies the algebraic Riccati equation $$A^{T}.K+K.A-K.B.R^{-1}.B^{T}.K=-Q$$
when $(q_{0},v_{0})$ satisfies $d_{g}^{2}(q_{0},q^{*})+\frac{\alpha}{k_{3}}|v_{0}|^{2}<inj_{q^{*}}(M)^{2}$. 
\end{theorem}
\end{minipage}}\\ \\
\textbf{Proof :}\\
We denote by $(q,v)\in TM$ the initial condition of the system, so the function\\ $V(q,v)=\inf_{u}J(u)$ is the solution of HJB equation $H(q,v,grad(V))=\gamma V$ with $$H=\inf_{u\in T_{q}M}<p_{1},v>+<p_{2},-grad(\Psi)+u>+U(q^{*},q)$$ $$+\frac{1}{2}|v|^{2}+\frac{\alpha}{2}|u-grad(\Psi)|^{2}$$
this function is convex in $u$ for all $(q,v)\in TM$ and $(p_{1},p_{2})\in T_{q}M^{2}$ in $T_{q}M$, she reach her minimum at $u=grad(\Psi)-\frac{1}{\alpha}p_{2}$, because $p_{2}=\nabla_{v}V$, solving the HJB equation will conclude.\\
So $H(q,v,grad(V))=\gamma V$ is equivalent to $$<\nabla_{q}V,v>-\frac{1}{2\alpha}|\nabla_{v}V|^{2}+U+\frac{1}{2}|v|^{2}=\gamma V.$$
A candidate function to solve HJB is $$V(q,v)=k_{1}U(q^{*},q)+\frac{k_{2}}{2}|v|^{2}+k_{3}<grad(U),v>,$$
we have $\frac{d}{dt}V(q,v)=<k_{1}grad(U),v>+<k_{2}v,\frac{Dv}{Dt}>+<k_{3}D_{v}grad(U),v>+<k_{3}grad(U),\frac{Dv}{Dt}>$ so $grad(V)=(k_{1}grad(U)+k_{3}D_{v}grad(U),k_{2}v+k_{3}grad(U))$ we replace $$\frac{1}{2}|grad(U)|^{2}=U,$$ $$D_{v}grad(U)=v,$$  to obtain $grad(V)=(k_{1}grad(U)+k_{3}v,k_{3}grad(U)+k_{2}v)$, putting this in the HJB equation gives $$k_{1}<grad(U),v>+k_{3}|v|^{2}-\frac{1}{2\alpha}[k_{3}^{2}|grad(U)|^{2}+k_{2}^{2}|v|^{2}$$ $$+2k_{3}k_{2}<grad(U),v>]+U+\frac{1}{2}|v|^{2}=\gamma V$$
$$(1-\frac{1}{\alpha}k_{3}^{2})U+\frac{1}{2}[1+2k_{3}-\frac{k_{2}^{2}}{\alpha}]|v|^{2}$$ $$+[k_{1}-\frac{k_{3}k_{2}}{\alpha}]<grad(U),v>=\gamma.V$$ so Riccati equations to solve for $k_{1,2,3}>0$ are $$1-\frac{k_{3}^{2}}{\alpha}=\gamma k_{1}$$ $$1+2k_{3}-\frac{k_{2}^{2}}{\alpha}=\gamma k_{2}$$ $$k_{1}-\frac{k_{3}k_{2}}{\alpha}=\gamma k_{3}$$
To write these equations as in the classical LQR theory, we put $K=\begin{pmatrix}k_{1} & k_{3} \\ k_{3} & k_{2} \end{pmatrix}$, $R=\alpha$, $Q=I_{2}$, et $B=\begin{pmatrix} 0 \\ 1 \end{pmatrix}$, $$A=\begin{pmatrix} -\gamma & 2\\ 0& -\gamma \end{pmatrix}.$$ Riccati equation takes the following form $$A^{T}.K+K.A-K.B.U^{-1}.B^{T}.K=-W,$$
that have unique solution that is positive definite. This gives $k_{1,2,3}>0$. It is clear that $C(A,B)=\begin{pmatrix}0& 2\\ 1& -\gamma \end{pmatrix}$ and so $(A,B)$ is controllable, and $Q,R>0$, so by the classical result in LQR theory, there exists unique solution positive definite for the Riccati equation, we recover the formula well known in LQR theory $$\begin{pmatrix}\frac{k_{3}}{\alpha}& \frac{k_{2}}{\alpha}\end{pmatrix}=R^{-1}B^{T}K.$$ On the other hand,  $$L(q,v)=\frac{k_{3}}{\alpha}U(q^{*},q)+\frac{1}{2}|v|^{2}$$ is a Lyapunov function of the optimal feedback. 
\subsection{Optimal Tracking}
Here we give an optimisation aspect of the tracking regulator exposed in [16] [8], we will show by computations on Riemannian manifolds and using the dynamic programming approach that the optimal regulator for a natural cost is the PD+FF regulator exposed by F. Bullo and R. Murray.\\
\fcolorbox{black}{lightgray}{
\begin{minipage}{0.48\textwidth}

\begin{theorem}{Optimisation Aspect of Tracking Regulator}\\
let $q_{ref}:[0,T]\rightarrow{M}$ be a smooth reference trajectory for the dynamic $\frac{Dq'}{Dt}=u$, consider the natural cost $$J(u)=\int_{0}^{T}[U(q_{ref},q_{u})+\frac{1}{2}|v_{u}-P(q_{u},q_{ref})v_{ref}|^{2}$$ $$+\frac{\alpha}{2}|u-u_{FF}|^{2}]dt.$$
So the optimal control is a feedback $u=u_{PD}+u_{FF}$ such that 
$$u_{FF}=\frac{d}{dt}P(q,q_{ref}(t))v_{ref}(t)+D_{v}P(q,q_{ref}(t))v_{ref}(t),$$
$$u_{PD}=-\frac{k_{3}(t)}{\alpha}\nabla U(q,q_{ref}(t))-\frac{k_{2}}{\alpha}(v-P(q,q_{ref}(t))v_{ref}(t)),$$
with $(k_{3}(t)/\alpha, k_{2}(t)/\alpha)=R^{-1}.B^{T}.K(t)$ where $K$ is the solution of the Riccati equation $$K'+A^{T}.K+K.A-K.B.R^{-1}.B^{T}.K=-Q$$ with the terminal condition $K(T)=0$.
\end{theorem}
\end{minipage}}\\ \\
\textbf{Proof :}\\
A candidate function to solve HJB is $$V(t,q,v)=k_{1}(t).U(q_{ref}(t),q)+\frac{k_{2}(t)}{2}|v-P(q,q_{ref})v_{ref}|^{2}$$ $$+k_{3}(t)<grad_{q}U(q_{ref},q),v-P(q,q_{ref})v_{ref}>.$$ Hamiltonian is  $$H(t,q,v,p)=\inf_{u\in T_{q}M}(<p,(v,u)>+U(q_{ref}(t),q)$$ $$+\frac{1}{2}|v-P(q,q_{ref}(t))v_{ref}(t)|^{2}+\frac{\alpha}{2}|u-u_{FF}(t,q,v)|^{2})$$ so $u=u_{FF}-\frac{1}{\alpha}p_{2}$, this gives $$H(t,q,v,p)=<p_{1},v>+<p_{2},u_{FF}>-\frac{1}{2\alpha}|p_{2}|^{2}$$ $$+U(q_{ref}(t),q)+\frac{1}{2}|v-P(q,q_{ref})v_{ref}|^{2}.$$
Some computations gives $$\frac{\partial V}{\partial t}=k_{1}'U+k_{2}'|v-P|^{2}+k_{3}'<\nabla U,v-P>$$ $$+k_{1}<\nabla_{1} U,v_{ref}>-k_{2}<\frac{d}{dt}P,v-P>$$ $$+k_{3}<\frac{d}{dt}\nabla_{2}U,v-P>+k_{3}<\nabla_{2}U,-\frac{d}{dt}P>$$
the compatibility of the pair $(U,P)$ gives $$<\nabla_{1}U,v_{ref}>=-<\nabla_{2}U,P>$$ 
\begin{figure}[H]
  \centering 
      \includegraphics[scale=0.7]{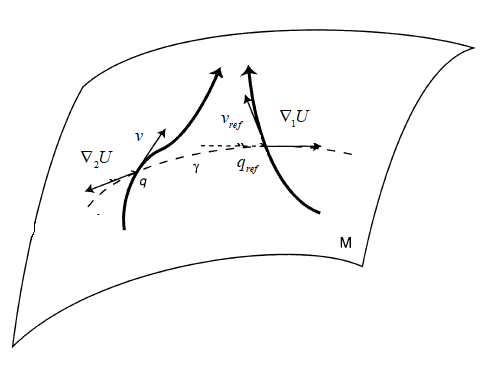} 
      \caption{Compatibility of the pair $(U,P)$ [22]}

           \label{fig:my_label}
            \end{figure}
now we compute $\nabla V$, for this $$\frac{d}{ds}V(t,q_{s},v_{s})=k_{1}\nabla_{2}U,v>+k_{2}<D_{s}v-D_{s}P,v-P>$$ $$+k_{3}<D_{s}\nabla_{2}U,v-P>+k_{3}<\nabla_{2}U,D_{s}v-D_{s}P>$$ en using the fact that$<P,P>$ is constant we conclude that $<D_{s}P,P>=0$, this gives $$\nabla_{q}V=k_{1}\nabla_{2}U-k_{2}D_{s}P+k_{3}(v-P)-I$$ $$\nabla_{v}V=k_{2}(v-P)+k_{3}\nabla_{2}U$$ such that $<I,v>=<\nabla U,D_{v}P>$, the Hamiltonian becomes 
$H=<k_{1}\nabla U-k_{2}D_{s}P+k_{3}(v-P)-I,v>+<k_{2}(v-P)+k_{3}\nabla U,u_{FF}>-\frac{1}{2\alpha}k_{2}^{2}|v-P|^{2}-\frac{1}{2\alpha}k_{3}^{2}|\nabla U|^{2}-\frac{k_{2}k_{3}}{\alpha}<\nabla U,v-P>+U+\frac{1}{2}|v-P|^{2}$ rearranging terms at HJB gives 
$$k_{1}'U+\frac{k_{2}'}{2}|v-P|^{2}+k_{3}'<\nabla U,v-P>$$ $$+k_{1}<\nabla U,v-P>-k_{2}<\frac{d}{d}P+D_{v}P,v-P>$$ $$+k_{3}<v-P,v>+k_{2}<u_{FF},v-P>+k_{3}<\nabla U,u_{FF}>$$ $$-k_{3}<\nabla U,\frac{d}{dt}P+D_{v}P>+k_{3}<\frac{d}{dt}grad_{q}U,v-P>$$ $$-\frac{1}{2\alpha}k_{2}^{2}|v-P|^{2}-\frac{1}{2\alpha}k_{3}^{2}|\nabla U|^{2}$$ $$-\frac{k_{2}k_{3}}{\alpha}<\nabla U,v-P>+U+\frac{1}{2}|v-P|^{2}=0$$  we conclude by taking $$u_{FF}(t,q,v)=\frac{d}{dt}P(q,q_{ref}(t))v_{ref}(t)+D_{v}P(q,q_{ref}(t))v_{ref}(t)$$ using the fact that $\frac{d}{dt}grad_{q}(U)=-P$, the equation becomes
$$k_{1}'U+\frac{k_{2}'}{2}|v-P|^{2}+k_{3}'<\nabla U,v-P>$$ $$+k_{1}<\nabla U,v-P>+k_{3}|v-P|^{2}-\frac{1}{2\alpha}k_{2}^{2}|v-P|^{2}$$ $$-\frac{1}{2\alpha}k_{3}^{2}|\nabla U|^{2}-\frac{k_{2}k_{3}}{\alpha}<\nabla U,v-P>+U+\frac{1}{2}|v-P|^{2}=0$$ this can be rewritten as
$$k_{1}'-\frac{k_{3}^{2}}{\alpha}+1=0,$$ 
$$k_{2}'+2k_{3}-\frac{k_{2}^{2}}{\alpha}+1=0,$$
$$k_{3}'+k_{1}-\frac{k_{2}k_{3}}{\alpha}=0,$$
or by the differential Riccati equation
$$K'(t)+A^{T}.K(t)+K(t).A-K(t).B.U^{-1}.B^{T}.K(t)=-W,$$
$$K(T)=0.$$
The optimal feedback is then 
  $$u(t,q,v)=-\frac{k_{3}(t)}{\alpha}\nabla U(q_{ref}(t),q)-\frac{k_{2}(t)}{\alpha}[v-P(q,q_{ref}(t))v_{ref}(t)]$$ with $(k_{3}/\alpha, k_{2}/\alpha)=R^{-1}B^{T}K(t)$.
\section{Applications and Simulations}
\subsection{Rigid body}
We apply the previous theoretical developments to a rigid body with fixed point, we recall that the rigid body is a very interesting mechanical system [21] [49]. The configuration space is the Lie group $SO(3)=\left\{R\in M_{3}(\mathbb{R}), R^{T}.R=R^{T}.R=I_{3}\right\}$.\\ 
The kinetic energy of the rigid body is left invariant, so we can compute the kinetic energy by his restriction $I$ to the tangent space $\mathfrak{so}(3)$ of $SO(3)$ at $I_{3}$, explicitly $\mathfrak{so}(3)=A_{3}(\mathbb{R})$. \\
By the canonical isomorphism $[.]^{\times}:(A_{3}(\mathbb{R}),[,])\rightarrow{(\mathbb{R}^{3},\times)}$ given by $$[\Omega]^{\times}=\frac{1}{2}\begin{bmatrix}
\Omega_{3,2}-\Omega_{2,3}\\ \Omega_{1,3}-\Omega_{3,1}\\ \Omega_{2,1}-\Omega_{1,2}
\end{bmatrix}$$
where $[,]$ is the commutator of matrices. \\
We denote by $[.]_{\times}$ the inverse of $[.]^{\times}$, in face these two isomorphisms preserves the Lie algebra structure, and then we can rewrite the dynamic equations 
$$R'=R\Omega$$ 
$$\omega'=I^{-1}(I.\omega\times\omega)+\tau$$
where $\omega=[\Omega]^{\times}$.
\begin{figure}[H]
  \centering 
      \includegraphics[scale=0.43]{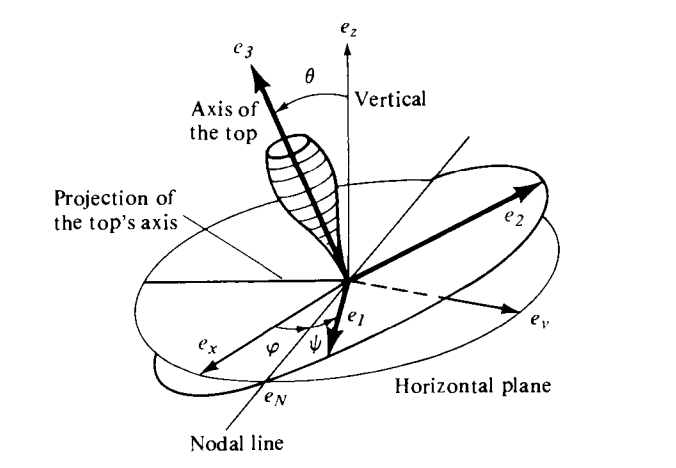} 
      \caption{Rigid body [8]}

           \label{fig:my_label}
            \end{figure}
We recall that 
$R=\alpha$, $Q=I_{2}$, et $B=\begin{pmatrix} 0 \\ 1 \end{pmatrix}$, $$A=\begin{pmatrix} -\gamma & 2\\ 0& -\gamma \end{pmatrix}$$\\
We simulate the theoretical developments using MATLAB, in order to solve the kinematic $R'=R\Omega$, we use the following Euler-Lie algorithm $$R_{i+1}=R_{n}\exp(h.\Omega_{i})$$
$$\omega_{i+1}=\omega_{i}+h.[I^{-1}(I.\omega_{i}\times\omega_{i})+\tau_{i}]$$
\subsection{Optimal regulation}
We are interested by the optimal regulation of the configuration to $R_{d}\in SO(3)$ by minimizing the cost $$J(\tau)=\int_{0}^{+\infty}\frac{1}{2}[d_{I}^{2}(R_{d},R)+|\Omega|_{I}^{2}+\alpha.|[\tau]_{\times}|_{I}^{2}].e^{-\gamma.t}dt$$
The optimal regulator is $$\tau=-k_{P}[\log(R_{d}^{T}R)]_{\times}-k_{D}\omega$$
for $(k_{P},k_{D})=lqr(A,B,Q,R)$.\\
 We take in the simulations $R_{d}=I_{3}$, $\gamma=-1$ and $\alpha=0.5$. This gives $(k_{P},k_{D})=(1.4142, 2.7671)$.\\ 
\\
\textbf{Simulations results :}\\ \\
Applying the previous control laws gives the following figures that shows the efficacy of the proposed optimal regulator, we see that the diagonal components converges to 1, and the other to 0, this conclude that $R\rightarrow{I}$. 
\begin{figure}[H]
  \centering 
      \includegraphics[scale=0.25]{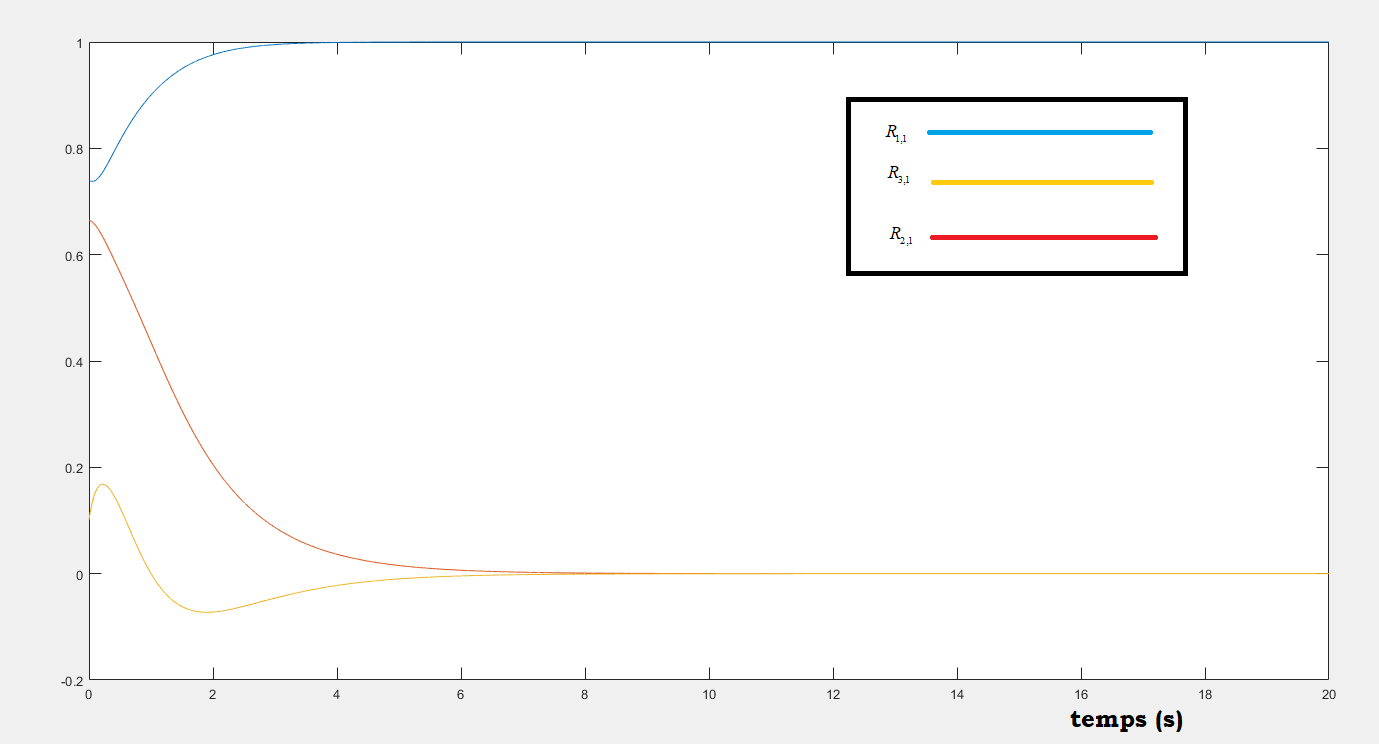} 
      \caption{Optimal regulation $R\rightarrow{I}$, 1st column}

           \label{fig:my_label}
            \end{figure}
            \begin{figure}[H]
  \centering 
      \includegraphics[scale=0.25]{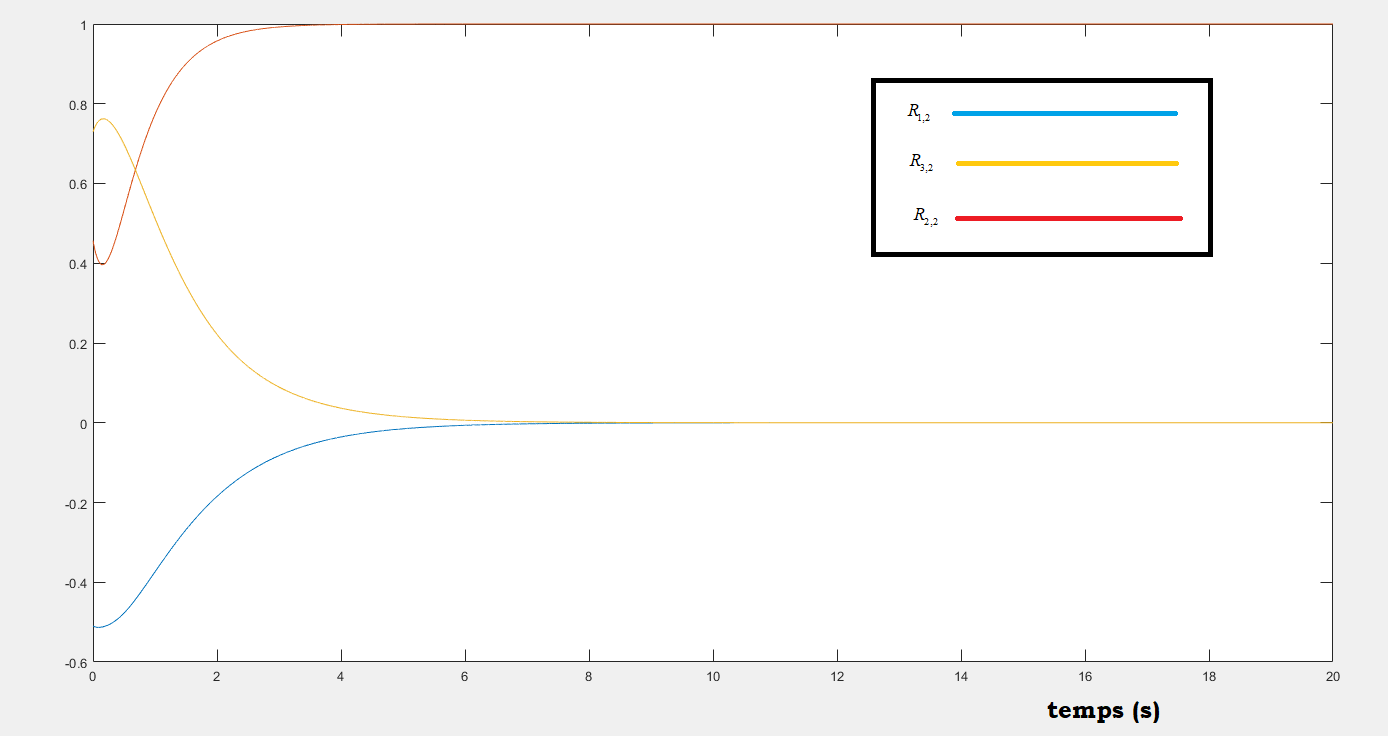} 
      \caption{Optimal regulation $R\rightarrow{I}$, 2nd column}

           \label{fig:my_label}
            \end{figure}
            \begin{figure}[H]
  \centering 
      \includegraphics[scale=0.25]{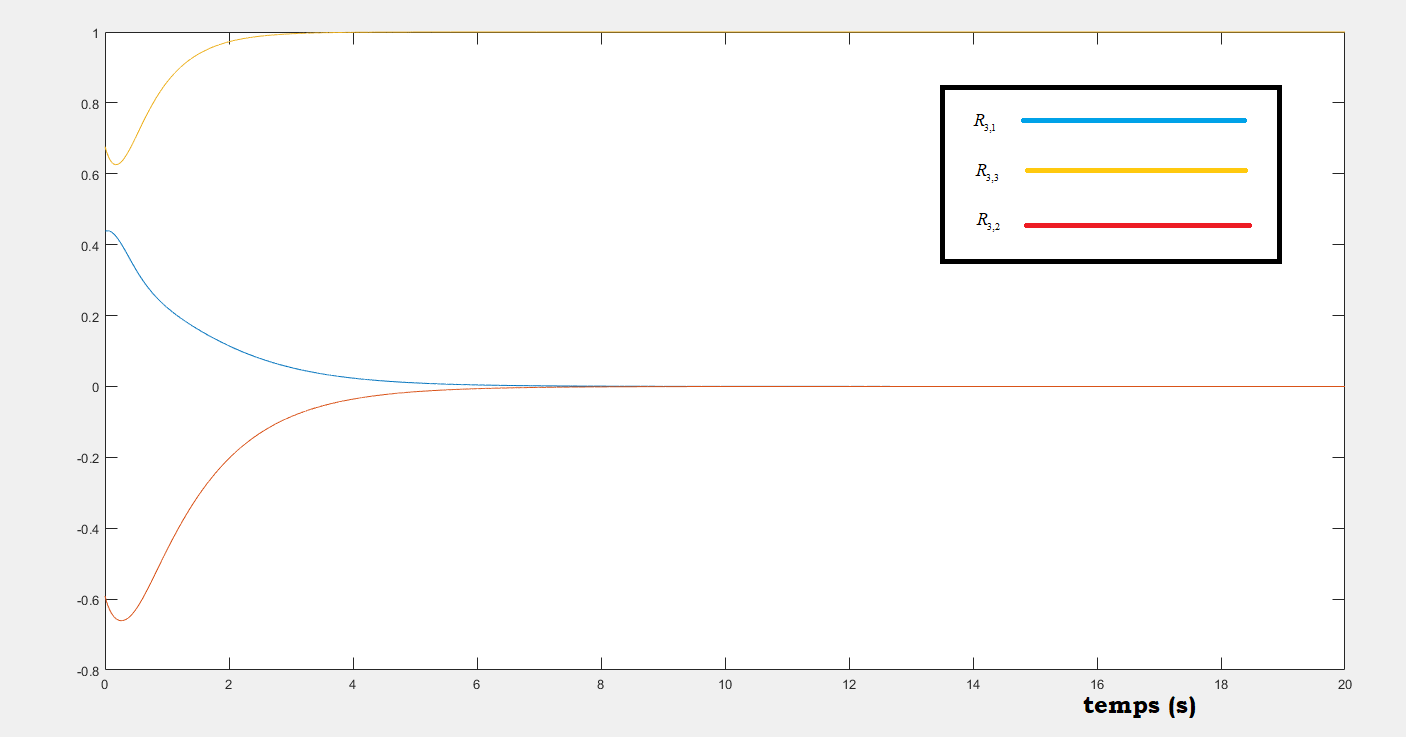} 
      \caption{Optimal regulation $R\rightarrow{I}$ 3rd column }

           \label{fig:my_label}
            \end{figure}
\subsection{Optimal tracking}
We are interested by the optimal tracking of the configuration to the smooth reference one $R_{ref}:[0,50]\rightarrow{SO(3)}$ by minimizing the cost 
$$J(\tau)=\int_{0}^{T}\frac{1}{2}[d_{I}^{2}(R_{ref},R)+|\Omega-P(R,R_{ref})(\Omega_{ref})|_{I}^{2}$$ $$+\alpha|[\tau-\tau_{FF}]_{\times}|_{I}^{2}]e^{-\gamma.t}dt.$$
When the metric is bi invariant, parallel transportation along the unique minimizing geodesic for two sufficiently close configurations is exactly the right translation, so we do this approximation 
$$P(R,R_{ref})\Omega_{ref}\cong \Omega_{ref}.R_{ref}^{T}.R$$
and by computations in [21] [26], the optimal regulator is
$$\tau_{FF}=\frac{1}{2}([\Omega]^{\times}\times R^{T}R_{ref}.[\Omega_{ref}]^{\times}-$$ $$I^{-1}.(I.R^{T}R_{ref}[\Omega_{ref}]^{\times}\times [\Omega]^{\times}+I.[\Omega]^{\times}\times R^{T}R_{ref}[\Omega_{ref}]^{\times})),$$
and the PD action is 
$$\tau_{PD}=-k_{P}[\log(R_{ref}^{T}R)]^{\times}-k_{D}.([\Omega]^{\times}-R^{T}R_{ref}.[\Omega_{ref}]^{\times}),$$
for $(k_{P},k_{D})=lqr(A,B,Q,R)$, we take $\gamma=-2$, $\alpha=1$. \\
So $(k_{P},k_{D})=(8.7852, 8.3357)$, we choose $\Omega_{ref}(t)=[0.5t,0.3t,0.4t]_{\times}$ and $R_{ref}$ is the unique trajectory such as $R_{ref}'=R_{ref}.\Omega_{ref}$. \\\\
\textbf{Simulations results :}\\ \\
The following figures shows the efficacy of the proposed regulator, we see that all the components of the attitude converge to the reference one. 

\begin{figure}[H]
  \centering 
      \includegraphics[scale=0.23]{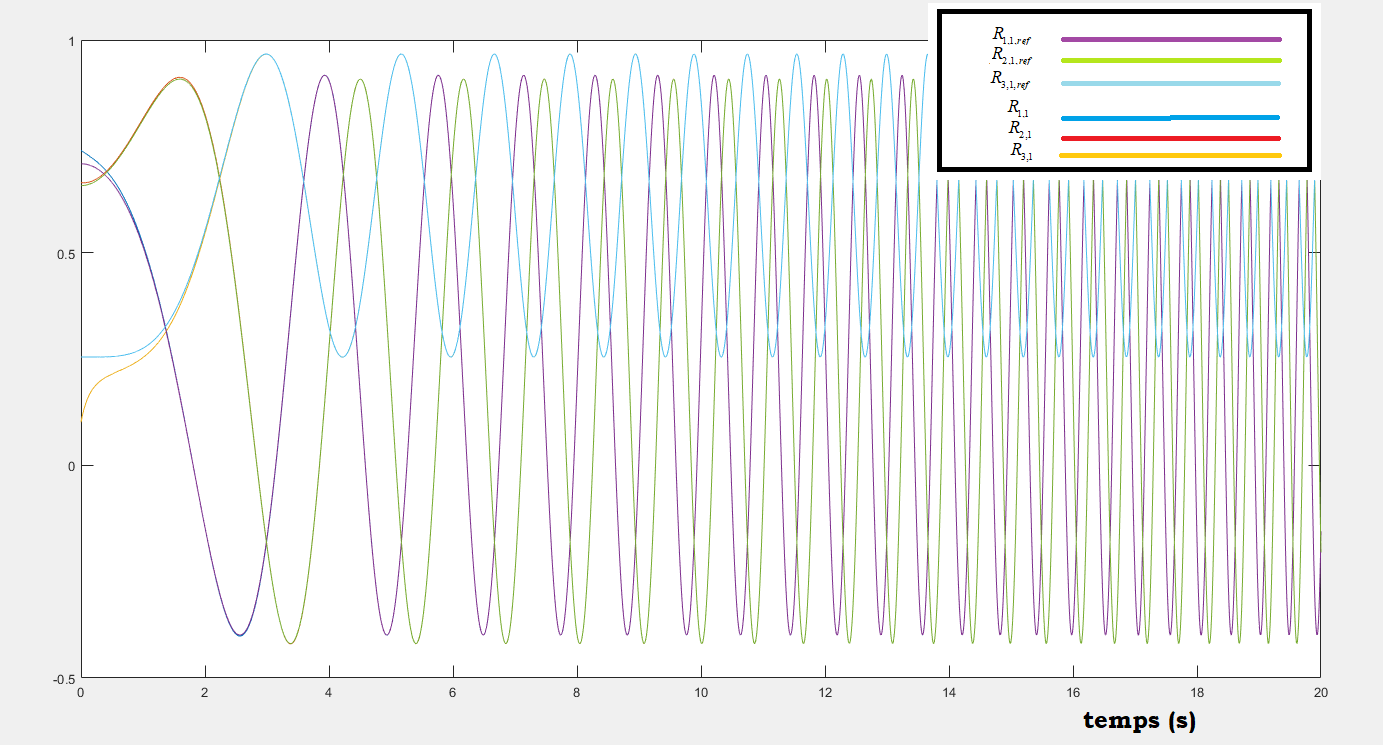} 
      \caption{Optimal tracking of $R_{ref}.[1,0,0]$}

           \label{fig:my_label}
            \end{figure}
            \begin{figure}[H]
  \centering 
      \includegraphics[scale=0.23]{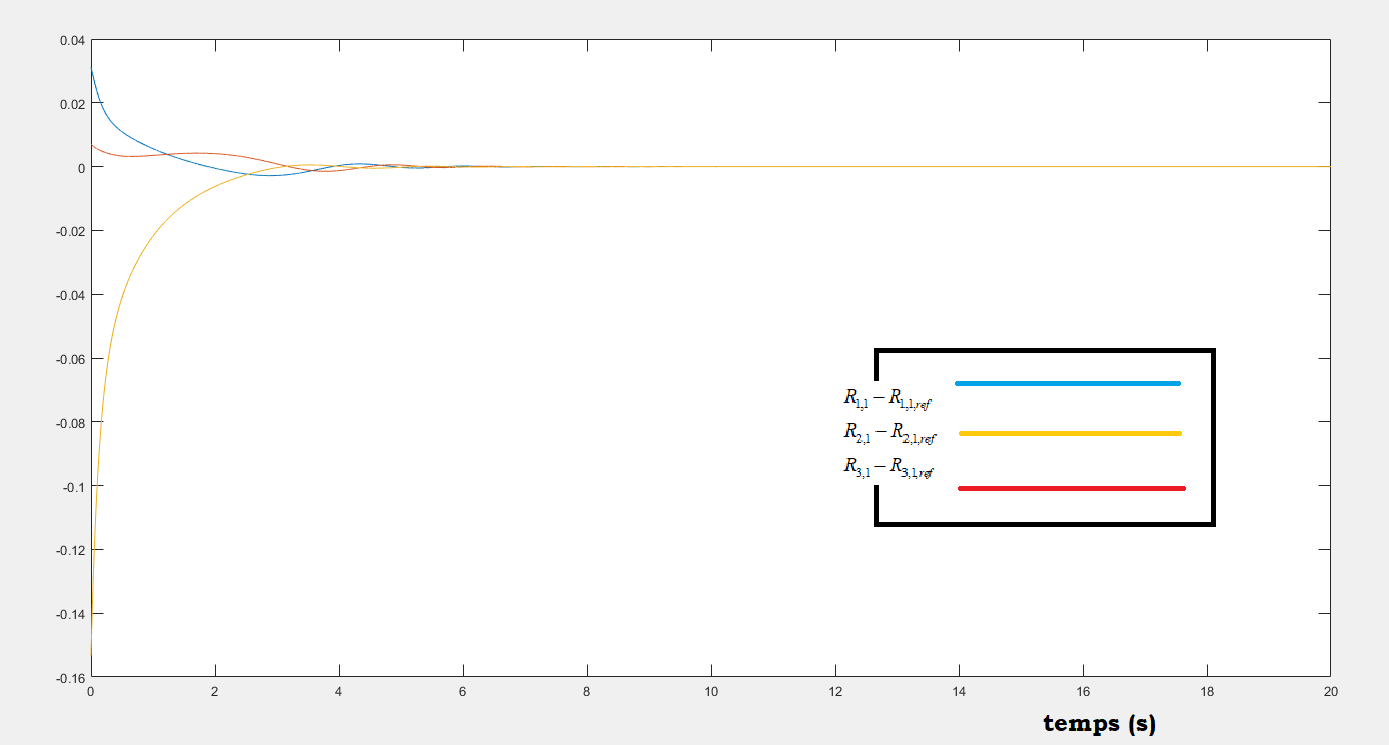} 
      \caption{Tracking error of $R_{ref}.[1,0,0]$}

           \label{fig:my_label}
            \end{figure}
\begin{figure}[H]
  \centering 
      \includegraphics[scale=0.23]{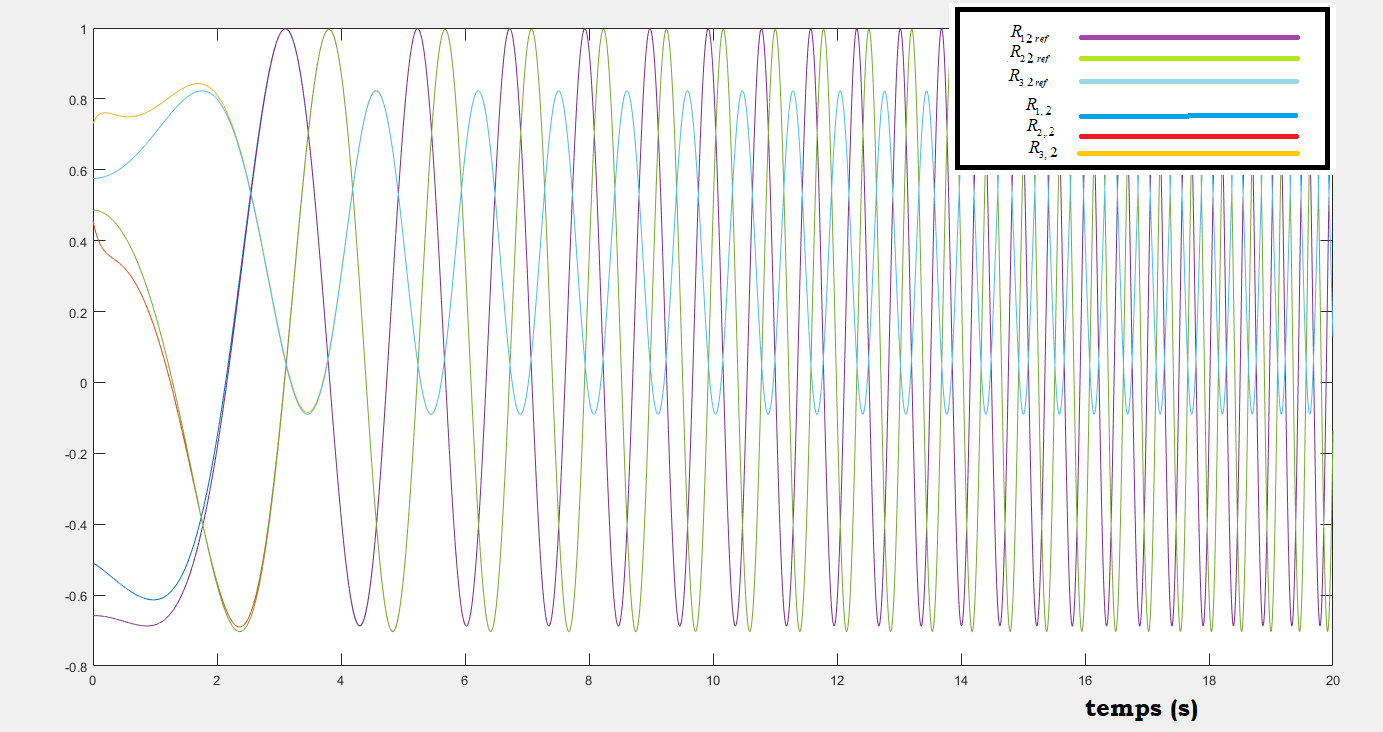} 
      \caption{Optimal tracking of $R_{ref}.[0,1,0]$}

           \label{fig:my_label}
            \end{figure}
            \begin{figure}[H]
  \centering 
      \includegraphics[scale=0.23]{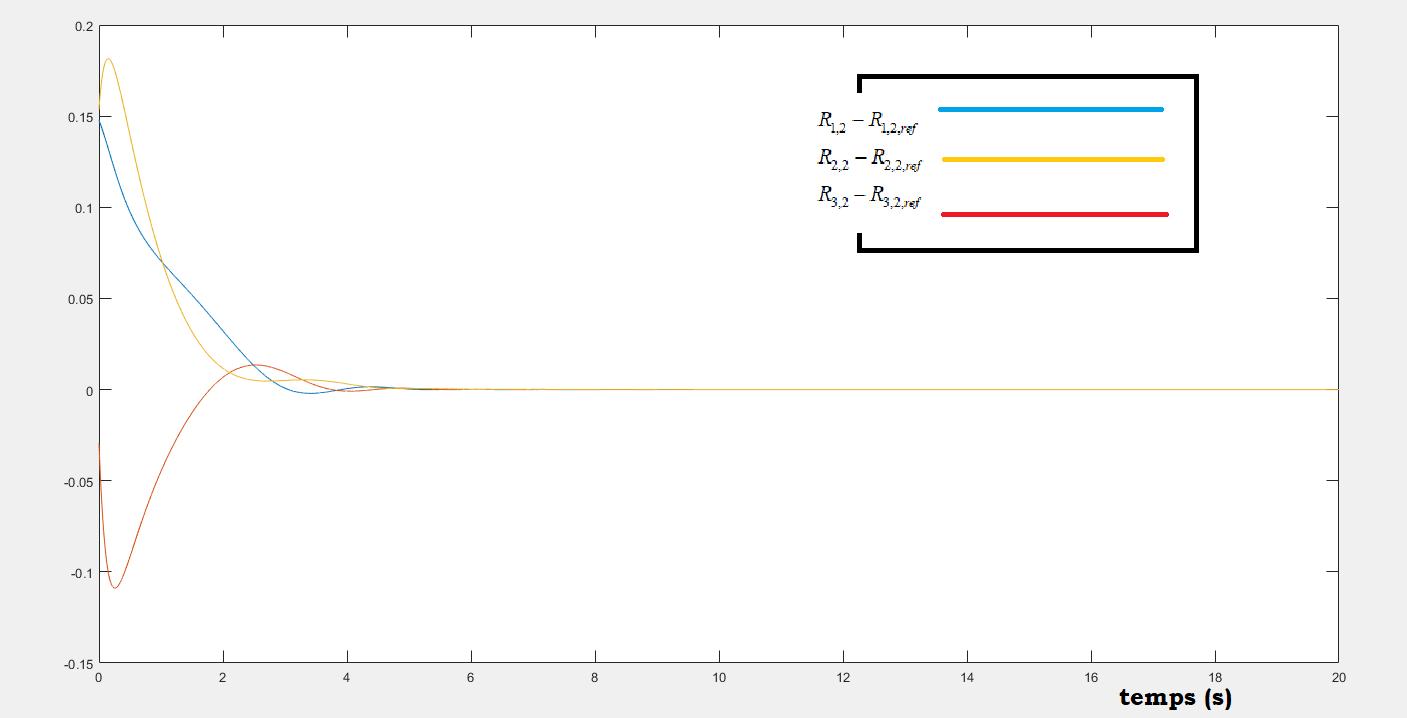}  
      \caption{Tracking error of $R_{ref}.[1,0,0]$}

           \label{fig:my_label}
            \end{figure}
            \begin{figure}[H]
  \centering 
      \includegraphics[scale=0.23]{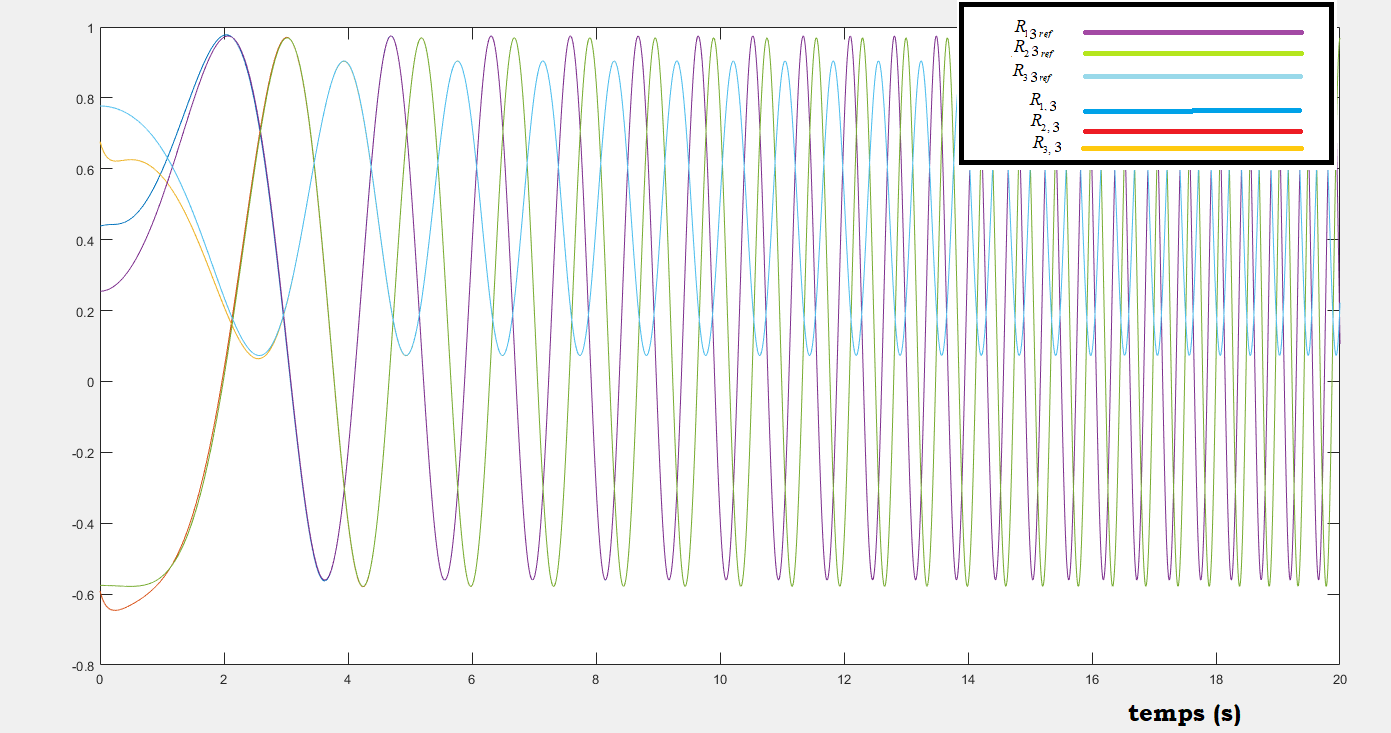} 
      \caption{Optimal tracking of $R_{ref}.[0,0,1]$}

           \label{fig:my_label}
            \end{figure}
            \begin{figure}[H]
  \centering 
      \includegraphics[scale=0.23]{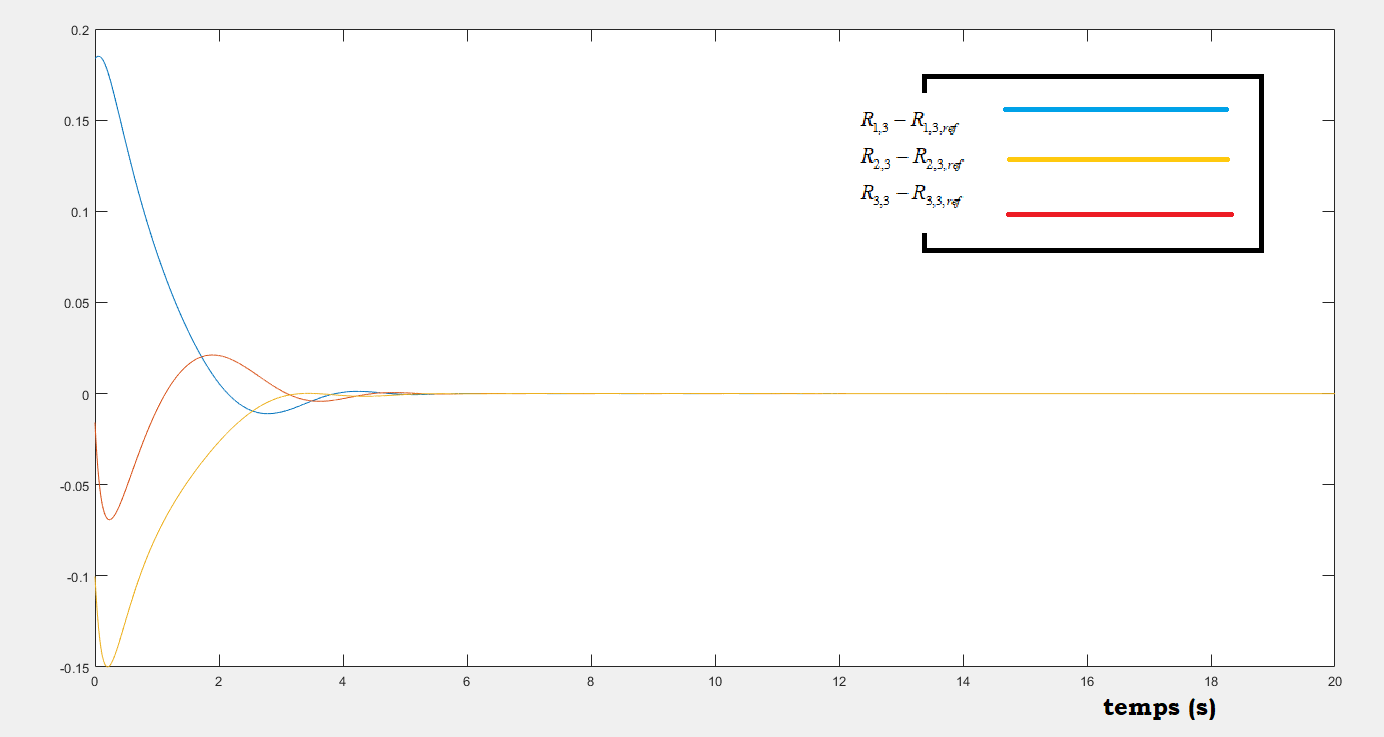} 
      \caption{Tracking error of $R_{ref}.[1,0,0]$}

           \label{fig:my_label}
            \end{figure}
\section{Conclusion}
Using functional analysis tools, we give some existence results about optimal control of robotic systems using an intrinsic formulation.\\ \\
We give an intrinsic formulation of PMP for robotic systems that involves explicitly the Riemannian curvature tensor, we then recover results of Bloch, Silva and Colombo [45] [46] [47]. \\
\\Finally the LQR extension for robotic systems with the geometric formulation is given using dynamic programming approach, and we give an optimisation aspect of the tracking regulation of F. Bullo and R. Murray in [21] [26].

\end{document}